\documentstyle[amstex,amssymb]{article}

\begin{document}

Assume that $\Theta $ is an arbitrary variety of groups. Let $W(X)$ be a
free group of the variety $\Theta $ over the finite set $X$ and $G$ is a
group in this variety ($G\in \Theta $). We can consider the ''affine space
over the group $G$'': ${\rm Hom}_\Theta (W(X),G)$. For every set $T\subset
W(X)$ we can consider the ''algebraic variety'' 
\[
A(T)_G=T_G^{\prime }=\{\mu \in {\rm Hom}_\Theta (W(X),G)\mid T\subset \ker
\mu \} 
\]
and $G$-closure of $T$ : 
\[
T_G^{\prime \prime }=\bigcap\limits_{\mu \in T_G^{\prime }}\ker \mu \subset
W(X). 
\]

{\bf Definition.} [Pl1]{\bf \ }{\it Groups }$G_1,G_2\in \Theta ${\it \ are
called }$X${\bf -equivalent}{\it \ if }$T_{G_1}^{\prime \prime
}=T_{G_2}^{\prime \prime }${\it \ holds for every set }$T\subset W(X)${\it .
Groups }$G_1,G_2\in \Theta ${\it \ are called {\bf geometrically equivalent}
(}denoted{\it \ $G_1\sim G_2$ ) if they are }$X${\it -equivalent for every
finite }$X${\it .}\vrule height7pt width4pt depth1pt\newline

It's clear that the relation \textquotedblright $\sim $\textquotedblright\
is an equivalence. It was proved that groups $G_{1},G_{2}\in \Theta $ are
geometrically equivalent if and only if every finitely generated subgroup of 
$G_{1}$ can be approximated by subgroup of $G_{2}$ and vice versa [Pl2]. We
denote{\bf \ }$G\prec H$ if the group $G$ can be approximated by the group $%
H $, i.e., there is $\{\varphi _{i}\mid i\in I\}\subset {\rm Hom}(G,H)$ such
that $\bigcap\limits_{i\in I}\ker \varphi _{i}=\left\{ 1\right\} $. It's
clear that the relation \textquotedblright $\prec $\textquotedblright\ is
the order. If we consider only nilpotent finitely generated groups, two
groups $G$ and $H$ are geometrically equivalent , iff $G\prec H$ and $H\prec
G$.

Also, it was proved that if two groups are geometrically equivalent, they
have the same identities [Pl1] and the same quasi-identities [PPT].\newline

{\bf Definitions.} {\it A variety }$\Theta ${\it \ of groups is called }{\bf %
Noetherian}{\it \ if the Noether chain condition holds for the normal
subgroups of every finitely generated free group }$W(X)${\it \ of this
variety.} [Pl3]

{\it A group }$G${\it \ in the arbitrary variety }$\Theta ${\it \ is called }%
[Pl3] {\bf geometrically Noetherian}{\it \ if for every finite set }$X${\it %
\ and every set }$T\subset W(X)${\it , where }$W(X)${\it \ is a free group
of the variety }$\Theta ${\it \ over the set }$X${\it , there is a finite
subset }$T_{0}\subset T${\it , such that }$T_{G}^{\prime }=\left(
T_{0}\right) _{G}^{\prime }${\it .}\vrule height7pt width4pt depth1pt\newline

For geometrically Noetherianity of $G$ we need the Noether chain condition
only for the $G$-closed normal subgroups in every finitely generated free
group $W(X)$ of $\Theta $, so every group $G$ in the Noetherian variety $%
\Theta $ is geometrically Noetherian.

The variety of a nilpotent class $s$ groups is Noetherian, because every
subgroup of nilpotent finitely generated group is finitely generated.By
[Pl3], two geometrically Noetherian groups are geometrically equivalent if
and only if they have the same quasiidentities. So two nilpotent class $s$
groups are geometrically equivalent if and only if they have the same
quasiidentities. The topic of quasiidentities and quasivarieties of
nilpotent groups (in most cases of nilpotent class $2$ groups) was
researched in many papers: [Is], [Fd1], [Fd2], [Bu], [Sh]. Now we can
resolve some questions in this topic by comparison of finitely generated
groups by the relation $\prec $. In particular, there is one-to-one
correspondence between classes of geometrically equivalent nilpotent groups
and quasivarieties generated by single nilpotent group.%

\section{Equivalence to Mal'tsev completion.}

It is known that for every nilpotent torsion free group $G$ there is a
Mal'tsev completion $\sqrt{G}$ - the minimal group, such that $G\subset 
\sqrt{G}$ and every $x\in \sqrt{G}$ has for every $n\in {\bf N}$ the $
x^{\frac 1n}\in \sqrt{G}$, such that $\left( x^{\frac 1n}\right) ^n=x$. The $
\sqrt{G}$ is the nilpotent group of the same class as $G$. The element $
x^{\frac 1n}\in \sqrt{G}$ is uniquely defined by $x\in \sqrt{G}$ and $n\in 
{\bf N}$.

Groups $G$ and $\sqrt{G}$ have the same identities ([Bau]). In the [PPT] the
question was asked when is the nilpotent torsion free group $G$
geometrically equivalent to its Mal'tsev completion $\sqrt{G}$ and, so,
groups $G$ and $\sqrt{G}$ have the same quasiidentities.\newline

{\bf Theorem 1.}

{\it If }$G${\it \ is a nilpotent class }$2${\it \ torsion free group, then
it is geometrically equivalent to its Mal'tsev completion }$\sqrt{G}${\it .}%
\newline

This theorem can be proved by simple calculation with the Mal'tsev basis and
Mal'tsev coordinats.\vrule height7pt width4pt depth1pt\newline

{\bf Corollary. }{\it A nilpotent class }$2$\ {\it torsion free group and
its Mal'tsev completion have same the quasiidentities.}\newline

This result was achieved independently also by Bludov and Gusev [BG]. In
these theses there is the example of nilpotent class $3$ torsion free group
(with $4$ generators) which is not geometrically equivalent to its Mal'tsev
completion. So this group and its Mal'tsev completion have different
quasi-identities. Therefore, the theorem of Baumslag cannot be extended to
quasi-identities.\newline

In {\bf Theorem 2}, we consider relatively free groups, i.e., groups which
are free in some subvariety of the variety of nilpotent class $s$ groups for
some $s\in {\bf N}$.\newline

{\bf Theorem 2. }{\it A nilpotent torsion free relatively free group is
geometrically equivalent to its Mal'tsev completion.}\newline

{\bf Corollary. }{\it A nilpotent torsion free relatively free group and its
Mal'tsev completion have same the quasiidentities.}\newline

In the proof of the {\bf Theorem 2} we used Lie ${\bf Q}$-algebras connected
with the nilpotent torsion free complete (in Mal'tsev sense) groups. All
these concepts we can see in [Ba]. First, in every nilpotent Lie ${\bf Q}$
-algebra $L$, we can define multiplication by Campbell-Hausdorff formula: 
\[
x\cdot y=x+y+\frac 12\left[ x,y\right] +\ldots 
\]
$L$ will be a complete nilpotent torsion free group by multiplication ''$
\cdot $'' (denoted as $L^{\circ }$). $L^n=\gamma _n\left( L^{\circ }\right) $
for every $n\in {\bf N}$, so the classes of nilpotency of $L$ and $L^{\circ
} $ coincide, and Abelian groups 
\[
\gamma _{n-1}\left( L^{\circ }\right) /\gamma _n\left( L^{\circ }\right)
\cong \left( L^{n-1}/L^n\right) ^{\circ }\cong L^{n-1}/L^n 
\]
are ${\bf Q}$-linear spaces. Conversely for every complete nilpotent torsion
free group $D$ there is nilpotent Lie ${\bf Q}$-algebra $L$, such that $
D\cong L^{\circ }$.

Secondly, if we have in some group the central filtration $\left\{ G_i\mid
i\in I\right\} $ ($I=\left\{ 1,\ldots ,n\right\} $ for some $n\in {\bf N}$,
or $I={\bf N}$), then we can consider the graded Lie ring $L$ defined by
this filtration: 
\[
\begin{array}{c}
L=\bigoplus\limits_{i\in I}G_i/G_{i+1}, \\ 
g_iG_{i+1}+h_iG_{i+1}=g_ih_iG_{i+1},\left[ g_iG_{i+1},g_jG_{j+1}\right]
=\left( g_i,g_j\right) G_{i+j+1}, \\ 
(g_i,h_i\in G_i,g_j\in G_j).%
\end{array}
\]
If Abelian groups $G_i/G_{i+1}$ are ${\bf Q}$-linear spaces, then $L$ is a
Lie ${\bf Q}$-algebra. If in the group $G$ there is central filtration $
\left\{ G_i\mid i\in I\right\} $ and $L\left( G\right) $ is the graded Lie
ring, defined by filtration $\left\{ G_i\mid i\in I\right\} $, and in the
group $H$ there is central filtration $\left\{ H_i\mid i\in I\right\} $ and $%
L\left( H\right) $ is the graded Lie ring, defined by filtration $\left\{
H_i\mid i\in I\right\} $, and $\varphi \in {\rm Hom}\left( G,H\right) $,
such that $G_i^\varphi \subset H_i$, then $\varphi $ induces $\overline{
\varphi }\in {\rm Hom}\left( L\left( G\right) ,L\left( H\right) \right) $
such that%
\[
\left( g_iG_{i+1}\right) ^{\overline{\varphi }}=g_i^\varphi G_{i+1}, 
\]
for every $g_i\in G_i$.

Every nilpotent class $s$\ torsion free group $G$ with $n$\ generators is
the factor-group of the $F_s\left( n\right) $ - the free nilpotent class $s$
\ group with $n$\ generators - by some normal isolated subgroup $T\lhd
F_s\left( n\right) ${\it . }Denote the natural homomorphism $\tau :F_s\left(
n\right) \rightarrow F_s\left( n\right) /T\cong G$.

For proving {\bf Theorem 2,} we must prove the following:\newline

{\bf Lemma.} {\it Let }$G${\it \ be a nilpotent class }$s${\it \ torsion
free group with }$n${\it \ generators, }$\sqrt{G}${\it \ its Mal'tsev
completion, }$H\leq \sqrt{G}${\it \ its finitely generated group. Then there
is }$k\in {\bf N}${\it , such that} 
\[
H\leq H_k=\left\langle \left( x_1^\tau \right) ^{\frac 1k},\ldots ,\left(
x_n^\tau \right) ^{\frac 1k}\right\rangle , 
\]
{\it \ where }$\left\{ x_1,\ldots ,x_n\right\} ${\it \ are free generators
of }$F_s\left( n\right) ${\it .}\newline

To prove this {\bf Lemma} we use the known theorem that every nilpotent
group $G$ can be generated by some subset $M\subset G$ and the commutant $%
\gamma _2\left( G\right) $ ($G=\left\langle M,\gamma _2\left( G\right)
\right\rangle $) can be generated by set $M$ ($G=\left\langle M\right\rangle 
$) (see [KM], 16.2.5).\newline

{\it Proof.}

Denote $L$ the Lie ${\bf Q}$-algebra, such that $L^{\circ }\cong \sqrt{G}$.
We can identify the elements of $L$, $L^{\circ }$ and $\sqrt{G}$. If $a,b\in 
\sqrt{G}=L$, $k\in {\bf N}$ , then, by the Campbell-Hausdorff formula, 
\[
a+b\equiv a\cdot b{\rm mod}\left[ L,L\right] , 
\]%
so 
\[
a^{\frac{1}{k}}\cdot b^{\frac{1}{k}}\equiv a^{\frac{1}{k}}+b^{\frac{1}{k}}=%
\frac{1}{k}\left( a+b\right) {\rm mod}\left[ L,L\right] , 
\]%
and 
\begin{equation}
a^{\frac{1}{k}}\cdot b^{\frac{1}{k}}\equiv \left( ab\right) ^{\frac{1}{k}}%
{\rm mod}\gamma _{2}\left( \sqrt{G}\right) .  \label{1}
\end{equation}

Every $g\in G$ has the form 
\[
g=\left( x_{i_{1}}^{\tau }\right) ^{\eta _{1}}\cdot \ldots \cdot \left(
x_{i_{r}}^{\tau }\right) ^{\eta _{r}}g_{2}, 
\]%
where $i_{1},\ldots ,i_{r}\in \left\{ 1,\ldots ,n\right\} $, $g_{2}\in
\gamma _{2}\left( G\right) <\gamma _{2}\left( \sqrt{G}\right) $, $\eta
_{1},\ldots ,\eta _{r}\in {\bf Z}$. Therefore, by (\ref{1}), 
\[
g^{\frac{1}{k}}\equiv \left( \left( x_{i_{1}}^{\tau }\right) ^{\eta
_{1}}\cdot \ldots \cdot \left( x_{i_{r}}^{\tau }\right) ^{\eta _{r}}\right)
^{\frac{1}{k}}\equiv \left( x_{i_{1}}^{\tau }\right) ^{\frac{\eta _{1}}{k}%
}\cdot \ldots \cdot \left( x_{i_{r}}^{\tau }\right) ^{\frac{\eta _{r}}{k}}%
{\rm mod}\gamma _{2}\left( \sqrt{G}\right) . 
\]%
Consequently, 
\[
\left\langle \left( x_{1}^{\tau }\right) ^{\frac{1}{k}},\ldots ,\left(
x_{n}^{\tau }\right) ^{\frac{1}{k}},\gamma _{2}\left( \sqrt{G}\right) \mid
k\in {\bf N}\right\rangle =\sqrt{G}. 
\]%
So 
\begin{equation}
\left\langle \left( x_{1}^{\tau }\right) ^{\frac{1}{k}},\ldots ,\left(
x_{n}^{\tau }\right) ^{\frac{1}{k}}\mid k\in {\bf N}\right\rangle =\sqrt{G}.
\label{2}
\end{equation}%
If $H=\left\langle h_{1},\ldots ,h_{m}\right\rangle $, where $h_{i}=\left(
x_{j_{1}}^{\tau }\right) ^{\frac{t_{i,1}}{k_{i,1}}}\cdot \ldots \cdot \left(
x_{j_{w_{i}}}^{\tau }\right) ^{\frac{t_{i,w_{i}}}{k_{i,w_{i}}}}$, $%
j_{1},\ldots ,j_{w_{i}}\in \left\{ 1,\ldots ,n\right\} $, $t_{i,l}\in {\bf Z}
$, $k_{i,l}\in {\bf N}$, $1\leq i\leq m$, then $H\leq H_{k}$, where $%
k=\prod\limits_{i,j}k_{i,j}$. The lemma is thus proved.\vrule height7pt
width4pt depth1pt\newline

{\it Proof of }{\bf Theorem 2. }

It is necessary to be proved only for a finitely generated group. Let $G$ be
a nilpotent class $s$ torsion free relatively free group with $n$
generators. $G=F_s\left( n\right) /T$, where $T$ is a isolated verbal
subgroup of $F_s\left( n\right) $. By the {\bf Lemma},{\bf \ }it is
necessary to approximate the%
\[
H_k=\left\langle \left( x_1^\tau \right) ^{\frac 1k},\ldots ,\left( x_n^\tau
\right) ^{\frac 1k}\right\rangle \leq \sqrt{G} 
\]
by $G$. The mapping 
\[
\varphi :x_i^\tau \rightarrow \left( x_i^\tau \right) ^k(1\leq i\leq n) 
\]
can be extended for the endomorphism of $G$ ( because $G$ is relatively free
) and $\sqrt{G}$ ([Ba], Chapter 6.6, Theorem 4). $H_k^\varphi \leq G$, so it
is only necessary to prove that $\ker \varphi =\left\{ 1\right\} $.

We shall consider for this purpose the graded Lie ${\bf Q}$-algebra $%
L_{\gamma }\left( \sqrt{G}\right) $ determined by central filtration $%
\left\{ \gamma _{i}\left( \sqrt{G}\right) \right\} _{i=0}^{s}$. The
endomorphism $\varphi \in {\rm End}\left( \sqrt{G}\right) $ induces
endomorphism $\overline{\varphi }\in {\rm End}\left( L_{\gamma }\left( \sqrt{%
G}\right) \right) $. By (\ref{2}), every 
\[
l_{i}\in \gamma _{i}\left( \sqrt{G}\right) /\gamma _{i+1}\left( \sqrt{G}%
\right) 
\]%
has the form: 
\[
l_{i}=\prod\limits_{j}w_{j}\cdot \gamma _{i+1}\left( \sqrt{G}\right) , 
\]%
where 
\[
w_{j}=w_{j}\left( \left( x_{r_{1}}^{\tau }\right) ^{q_{j,1}},\ldots ,\left(
x_{r_{i}}^{\tau }\right) ^{q_{j,i}}\right) 
\]%
( $r_{1},\ldots ,r_{i}\in \left\{ 1,\ldots ,n\right\} $, $q_{j,l}\in {\bf Q}$
) is a group commutator with length $i$. So 
\[
l_{i}^{\overline{\varphi }}=\sum\limits_{j}\left( w_{j}^{\prime }\right) ^{%
\overline{\varphi }}=\sum\limits_{j}k^{i}w_{j}^{\prime }=k^{i}l_{i} 
\]%
and $\ker \overline{\varphi }=\left\{ 0\right\} $, because $\gamma
_{i}\left( \sqrt{G}\right) /\gamma _{i+1}\left( \sqrt{G}\right) $ are ${\bf Q%
}$-linear spaces. On the other hand, if $g_{i}\in \gamma _{i}\left( \sqrt{G}%
\right) \cap \ker \varphi $ ( $0\leq i\leq s$ ), then $g_{i}\gamma
_{i+1}\left( \sqrt{G}\right) \in \ker \overline{\varphi }$. Therefore $\ker
\varphi =\left\{ 1\right\} $. The proof is complete.\vrule height7pt
width4pt depth1pt\newline


\section{Groups with the small rank of center.}


It is a very interesting problem to classify the classes of geometrically
equivalent groups in the variety of nilpotent groups of some fixed class $s$%
, or, in other words, to describe all quasivarietes generated by single
group of this variety.

It was proved by A. Berzins [Be], that two Abelian groups are geometrically
equivalent if and only if for every prime number $p$ the exponents of their
corresponding $p$-Sylow subgroups coincide, and if one of these group is not
periodic, then second group is not periodic either. The second step in this
program will be the classifying of classes of geometrically equivalent
nilpotent torsion free class $2$ groups. This step can be considered as an
approach to resolving the very sophisticated problem of classifying of all
nilpotent torsion free class $2$ groups.

In the researching of the geometrical equivalence of nilpotent torsion free
class $2$ groups we can (by {\bf Theorem 1}) consider only complete
nilpotent class $2$ groups, which are completions of finitely generated
torsion free groups. We shall say that these complete group have a finite
rank.

We used the approach of [GrS], which permits us to consider the problem of
approximating complete nilpotent class $2$ torsion free groups of finite
rank by the technique of linear algebra. Let $A_{0}$ be the nilpotent class $%
2$ finitely generated torsion free group and $A$ its Mal'tsev completion. We
say, that the group $A$ is the nilpotent class $2$ torsion free complete
group of finite ranks. We can identify the the nilpotent class $2$ torsion
free complete group $A$ of finite ranks with the graded Lie ${\bf Q}$%
-algebra 
\[
L_{\zeta }\left( A\right) =A/Z\left( A\right) \oplus Z\left( A\right) , 
\]%
or with the the pair of ${\bf Q}$ - vector space 
\[
\begin{array}{c}
V_{A}=A/Z\left( A\right) (\dim V_{A}={\rm rank}A_{0}/Z\left( A_{0}\right)
<\infty ), \\ 
W_{A}=Z\left( A\right) (\dim W_{A}={\rm rank}Z\left( A_{0}\right) <\infty )-%
\end{array}%
\]%
and the nonsingular alternating bilinear mapping 
\[
\begin{array}{c}
\lbrack \_,\_]_{A}:V_{A}\times V_{A}\ni (v_{1}Z\left( A\right) ,v_{2}Z\left(
A\right) )\rightarrow \\ 
\rightarrow \lbrack v_{1}Z\left( A\right) ,v_{2}Z\left( A\right)
]_{A}=\left( v_{1},v_{2}\right) \in Z\left( A\right) =W_{A},%
\end{array}%
\]%
which define Lie brackets in the $L_{\zeta }\left( A\right) $.

If holds $A\sim A_1$, $B\sim B_1$ ($A$, $A_1$, $B$, $B_1$ - groups), then $%
A\times B\sim A_1\times B_1$. So in this study we can consider only groups
nondecomposable in the direct product. So, we can assume for a complete
nilpotent class $2$\ group $A$ that $Z\left( A\right) =\left[ A,A\right] $,
or 
\begin{equation}  \label{3}
{\rm im}[\_,\_]_A=W_A.
\end{equation}
Otherwise $A$ is decomposable in the direct product.

Let $A=(V_{A},W_{A},[\_,\_]_{A})$ and $B=(V_{B},W_{B},[\_,\_]_{B})$ be two
complete nilpotent class $2$ groups, which fulfill condition (\ref{3}). The
homomorphism $\Phi :A\rightarrow B$ we can consider as the pair $\left(
\varphi ,\psi \right) $ of two linear maps

\[
\begin{array}{c}
\varphi :V_{A}=A/Z\left( A\right) \ni aZ\left( A\right) \rightarrow a^{\Phi
}Z\left( A\right) \in B/Z\left( B\right) =V_{B} \\ 
\psi :W_{A}=Z\left( A\right) \ni z\rightarrow z^{\Phi }\in Z\left( B\right)
=W_{B},%
\end{array}
\]
which makes the following diagram commutative 
\[
\begin{array}{ccc}
\Lambda ^{2}V_{A} & \overset{[\_,\_]_{A}}{\rightarrow } & W_{A} \\ 
\Lambda ^{2}\varphi \downarrow &  & \downarrow \psi \\ 
\Lambda ^{2}V_{B} & \overset{[\_,\_]_{B}}{\rightarrow } & W_{B}%
\end{array}%
, 
\]
i.e., they satisfy the condition $\psi \lbrack \_,\_]_{A}=[\varphi
\_,\varphi \_]_{B}$.

It is clear that if $A=(V_{A},W_{A},[\_,\_]_{A})$ and $B=(V_{B},W_{B},[\_,%
\_]_{B})$ are two complete nilpotent class $2$ groups of finite ranks, then $%
A\prec B$ there is the famili of pairs of linear maps $\{(\varphi _{i},\psi
_{i})\mid i\in I\}$ such that

1)$\varphi _i\in {\rm Hom}_{{\bf Q}}(V_A,V_B)$,

2)$\psi _i\in {\rm Hom}_{{\bf Q}}(W_A,W_B)$,

3)$\psi _i([\_,\_]_A)=[\varphi _i\_,\varphi _i\_]_B$,

4)$\bigcap\limits_{i\in I}\ker \psi _i=0$.

Denote $S(k,V)$ ( $2\mid k$, $\dim V=s\geq k$ ) the alternate bilinear form: 
$\Lambda ^2V\rightarrow {\bf Q}$ such that, if $\{e_1,...,e_n\}$ is the
basis of $V$, then 
\[
S(k,V)(e_{2i-1},e_{2i})=-S(k,V)(e_{2i},e_{2i-1})=1(1\leq i\leq \frac k2), 
\]
otherwise 
\[
S(k,n)(e_i,e_j)=0. 
\]
Denote the group $(V,W,S(k,V))$, where $\dim V=k$, $2\mid k$, $\dim W=1$ as $%
N(k,1)$. It is clear that all Mal'tsev completions of nilpotent class $2$
finitely generated torsion free group with the cyclic center are isomorphic
to the group $N(k,1)$ ( $2\mid k$ ).

Our aim is to prove the following:\newline

{\bf Theorem 3. }{\it Let }$G_{1}${\it , }$G_{2}${\it \ two nilpotent
torsion free class }$2${\it \ finitely generated groups with the cyclic
center. Then }$G_{1}${\it \ and }$G_{2}$ {\it are geometrically equivalent (}%
$G_{1}\sim G_{2}${\it ) if and only if their Mal'tsev completions are
isomorphic (}$\sqrt{G_{1}}\cong \sqrt{G_{2}}${\it ).}

and the\newline

{\bf Theorem 4. }{\it Let }$G_{1}${\it , }$G_{2}${\it \ two nilpotent
torsion free class }$2${\it \ finitely generated groups, whose centers have
rank }$2${\it . Then }$G_{1}\sim G_{2}${\it \ if and only if, or there is
nilpotent torsion free class }$2${\it \ finitely generated group with the
cyclic center }$N${\it , such that }$G_{1}\sim N\sim G_{2}${\it , or }$\sqrt{%
G_{1}}\cong \sqrt{G_{2}}${\it .}

By {\bf Theorem 3,} we have that there is a countable subset of
quasi-varaietes in the quasi-variety of the nilpotent class $2$ torsion free
groups. This result was achieved in 1975 by Isakov [Is] by a complicated
calculation, and here we achieve it easily.\newline

{\bf Proposition 1. }{\it Let }$A=(V_A,W_A,[\_,\_]_A)$,{\it \ the complete
nilpotent class }$2${\it \ group of the finite rank. Then }$N(k,1)\prec A$%
{\it \ iff there is }$V\leq V_A${\it \ such that }$\dim V=k${\it , }$\dim
[V,V]_A=1${\it \ and the bilinear form }$[\_,\_]_{A\mid V\times V}${\it \ is
nonsingular.}\newline

{\it Proof.}

$N(k,1)\prec A$ iff the $N(k,1)$ is embedded in the $A$.\vrule height7pt
width4pt depth1pt\newline

{\bf Corollary.} {\it If }$k_{1}\neq k_{2}${\it \ ( }$2\mid k_{1},k_{2}${\it %
\ ), then }$N(k_{1},1)${\it \ is not geometrically equivalent to the }$%
N(k_{2},1)${\it .}\vrule height7pt width4pt depth1pt\newline

{\bf Theorem 3} is proved.\vrule height7pt width4pt depth1pt\newline

{\bf Proposition 2. }{\it Let }$A=(V_A,W_A)${\it \ the complete nilpotent
class }$2${\it \ group of the finite rank. Then }$A\prec N(k,1)${\it \ iff
there are }$r=\dim W_A${\it \ linear independent functionals }$\psi _1$ {\it %
,...,}$\psi _r\in {\rm Hom}_{{\bf Q}}(W_A,{\bf Q)}${\it , such that the } $%
{\rm rank}\psi _i([\_,\_]_A)\leq k${\it \ for every }$i\in \{1,...,r\}$ {\it %
.}\newline

{\it Proof.}

Functionals{\it \ }$\psi _{1},...,\psi _{r}\in {\rm Hom}_{{\bf Q}}(W_{A},%
{\bf Q)}$ are linear independent iff $\bigcap\limits_{i=1}^{r}\ker \psi
_{i}=0$.\vrule height7pt width4pt depth1pt\newline

Let $A$, $B$ two complete nilpotent class $2$ groups of finite ranks. We say
that the comparison $A\prec B$ is realized by functionals if there is $%
\{\Phi _i=(\varphi _i,\psi _i)\mid i\in I\}\subset {\rm Hom}(A,B)$, such
that $\dim {\rm im}\psi _i=1$ for every( $i\in I$ ), $\bigcap\limits_{i\in
I}\ker \psi _i=0$.\newline

{\bf Corollary.} {\it Let }$A=(V_{A},W_{A},[\_,\_]_{A})${\it \ and } $%
B=(V_{B},W_{B},[\_,\_]_{B})${\it \ be two complete nilpotent class 2 groups
of finite ranks and the comparison }$A\prec B${\it \ realized by
functionals. Then there is }$k\in {\bf N}${\it , such that }$A\prec
N(k,1)\prec B${\it .}\newline

{\it Proof. }

Let $k=\max \limits_{i\in I}$ ${\rm rank}\psi _i[\_,\_]_A$. By {\bf %
Proposition 2}, we have $A\prec N(k,1)$. Without loss of generality,%
\[
{\rm rank}\psi _1([\_,\_]_A)=k\leq \dim V_A=a. 
\]
We have $[\varphi _1\_,\varphi _1\_]_B=\psi _1([\_,\_]_A)$ and $%
[V_A,V_A]_A=W_A$, so 
\[
\begin{array}{c}
\dim [\varphi _1(V_A),\varphi _1(V_A)]_B=\dim \psi _1([V_A,V_A]_A)=1, \\ 
{\rm rank}[\_,\_]_{B\mid \varphi _1(V_A)\times \varphi _1(V_A)}=k.%
\end{array}
\]
Thus there is the vector space $V\leq \varphi _1(V_A)\leq V_B$ such that 
\[
\dim V=k, \dim [V,V]_B=1 
\]
and the bilinear form $[\_,\_]_{B\mid V\times V}$ is nonsingular. Therefore,
by {\bf Proposition 1}, we have $N(k,1)\prec B$.\vrule height7pt width4pt
depth1pt\newline

{\bf Proposition 3.} {\it Let }$A=(V_A,W_A,[\_,\_]_A)${\it , } $%
B=(V_B,W_B,[\_,\_]_B)${\it \ two complete nilpotent class 2 torsion free
groups of finite ranks, }${\rm rank}Z(A)=2${\it . Then\ }$A\prec B${\it \
iff }

{\it 1)there is }$k\in {\bf N}${\it , such that }$A\prec N(k,1)\prec B$ {\it %
, or}

{\it 2)there is }$\Phi =(\varphi ,\psi )\in {\rm Hom}(A,B)${\it , such that }
$\ker \psi =0${\it , }$\ker \varphi =0${\it .}

{\it Proof. }

The condition 2) is equal to the embedding of the group $A$ into group $B$,
so condition 1), as well as condition 2) imply $A\prec B$.

The comparison $A\prec B$ can be realized by functionals or by embedding
(because ${\rm rank}Z(A)=2$). So, by {\bf Corollary} from {\bf Proposition 2,%
} this comparison implies condition 1), or condition 2).\vrule height7pt
width4pt depth1pt\newline

{\bf Proposition 4.}{\it \ Let }$A=(V_{A},W_{A},[\_,\_]_{A})${\it \ and }$%
B=(V_{B},W_{B},[\_,\_]_{B})${\it \ be two complete nilpotent class }$2${\it %
\ torsion free groups of finite ranks, whose centers have rank }$2${\it . If 
}$A\sim B${\it , then or there is }$k\in 2{\bf N}${\it \ such that }$A\sim
B\sim N(k,1)${\it , or }$A\cong B${\it .}\newline

{\it Proof.}

If $A\sim B$, then $A\prec B$, so, by {\bf Proposition 3}, or there is $k\in
2{\bf N}$, such that $A\prec N(k,1)\prec B$, or there is $\Phi =(\varphi
,\psi )\in \hom (A,B)$, such that $\ker \varphi _{A}=0$, $\ker \psi _{A}=0$.
In the first case, we have $A\sim B\sim N(k,1)$, because $A\sim B$. In the
second case, the proof is complete by symmetry.\vrule height7pt width4pt
depth1pt\newline

The {\bf Theorem 4} is the corollary of the {\bf Proposition 4}.

{\bf Proposition 5. }{\it Let }$A=(V_{A},W_{A})${\it \ the complete
nilpotent class }$2${\it \ group of the finite rank, which center has the
rank }$2$ {\it . Let }$\dim V_{A}=a${\it , }$2\mid a${\it . It cannot be
that }$A\sim N\left( a,1\right) ${\it .}\newline

{\it Proof. }

It follows from (\ref{3}) and {\bf Proposition 1}.\vrule height7pt width4pt
depth1pt\newline

In [Is] the question was posed: Are there other quasi-varieties of nilpotent
class $2$ groups besides the quasi-varieties generated by groups $N(k,1)$?
In [Bu] it was proved by complicated calculation that there are a continuum
of quasi-varieties of nilpotent class $2$ torsion free groups. By our method
we can easily answer the question of [Is] as ''yes''.\newline

{\bf Example.}

The group defined by 
\[
\left[ \_,\_\right] _A=\left( \left( 
\begin{array}{cccc}
0 & 1 & 0 & 0 \\ 
-1 & 0 & 0 & 0 \\ 
0 & 0 & 0 & 1 \\ 
0 & 0 & -1 & 0%
\end{array}
\right) ,\left( 
\begin{array}{cccc}
0 & 0 & 1 & 0 \\ 
0 & 0 & 0 & 0 \\ 
-1 & 0 & 0 & 0 \\ 
0 & 0 & 0 & 0%
\end{array}
\right) \right) , 
\]
i. e., the basis of $V_A$ is $\left\{ v_1,v_2,v_3,v_4\right\} $ and the
basis of $W_A$ is $\left\{ w_1,w_2\right\} $, $\left[ v_1,v_2\right] _A=%
\left[ v_3,v_4\right] _A=w_1$, $\left[ v_1,v_3\right] _A=w_2$, is not
geometrically equivalent to any group $N\left( k,1\right) $.

Actually, it is clear that $A$ can be equivalent only to $N\left( 2,1\right) 
$ or $N\left( 4,1\right) $. $A$ is not equivalent to the $N\left( 4,1\right) 
$ by {\bf Proposition 5}.{\bf \ }If $A\sim N\left( 2,1\right) $, then by 
{\bf Proposition 2} there must be two linear independent functionals $\psi
_{1},\psi _{2}\in W_{A}^{\ast }$, such that ${\rm rank}\psi
_{i}([\_,\_]_{A})\leq 2<4$ ($i=1,2$). Let $\left\{ \chi _{1},\chi
_{2}\right\} $ is the basis of $W_{A}^{\ast }$ dual for the $\left\{
w_{1},w_{2}\right\} $. Let 
\[
\psi =\lambda _{1}\chi _{1}+\lambda _{2}\chi _{2}\in W_{A}^{\ast }, 
\]%
i.e., 
\[
\begin{array}{c}
\psi ([\_,\_]_{A})=\left( 
\begin{array}{cccc}
0 & \lambda _{1} & \lambda _{2} & 0 \\ 
-\lambda _{1} & 0 & 0 & 0 \\ 
-\lambda _{2} & 0 & 0 & \lambda _{1} \\ 
0 & 0 & -\lambda _{1} & 0%
\end{array}%
\right) . \\ 
\det \left( 
\begin{array}{cccc}
0 & \lambda _{1} & \lambda _{2} & 0 \\ 
-\lambda _{1} & 0 & 0 & 0 \\ 
-\lambda _{2} & 0 & 0 & \lambda _{1} \\ 
0 & 0 & -\lambda _{1} & 0%
\end{array}%
\right) =\lambda _{1}^{4}=0\Longleftrightarrow \lambda _{1}=0,%
\end{array}%
\]%
so there is only one $\psi \in W_{A}^{\ast }$ such that ${\rm rank}\psi
_{i}([\_,\_]_{A})<4$. $A$ is not equivalent to the $N\left( 2,1\right) $.

I would like to thank B. Plotkin for the attention towards my researches in
this topics.

{\bf References:}

[Ba] Bahturrin Ju.A. {\it Lectures on Lie Algebras,} Academie-Verlag,
Berlin, 1978.

[Bau] Baumslag G. {\it On the residual nilpotence of some varietal product,}
Trans. Amer. Math. Soc. {\bf 109}, 1963, 357-365

[Be] Berzins A. {\it Geometrical equivalence of algebras,} International
Journal of Algebra and Computations (1998).

[BG] Bludov V.V., Gusev B.V., {\it About geometrical equivalence of groups,}
Theses of report.

[Bu] Budkin A.I., {\it A lattice of quasivarieties of nilpotent groups,}
Algebra and Logic, {\bf 33} (1), 1994.

[Fd1] Fedorov A.N., {\it O podkvzimnogoobrazijah nil'potentnyh minimal'nyh
neabelevyh mnogoobrazij grupp,} Sib. Mat. Zhurnal, {\bf 21} (6), 1980.

[Fd1] Fedorov A.N., {\it Kvazitozhdestva svobodnoj }$2${\it -nilpotentnoj
gruppy,} Mat. Zap., {\bf 40} (5), 1986

[GrS] Grunewald F.J., Scharlau R. {\it A note on finitely generated torsion
free nilpotent groups of class }$2${\it ,} Journal of Algebra. {\bf 58,} 1979

[Is] Isakov R.A. {\it O kvzimnogoobrazijah }$2${\it -stepenno nilpotentnyh
grup,} Izv. Aka\-demii Nauk UzSSR, 1976, 3, Serija fiz-mat. nauk.

[KM] Kargaplov M.I., Merzljakov Ju.I. {\it Fundamentals of the Theory of
Groups,} Springer Verlag, New York, 1979.

[Pl1] Plotkin B. {\it Algebraic logic, varieties of algebras and algebraic
varieties,} Proc. Int. Alg. Conf., St. Petersburg, 1995, Walter Gruyter, New
York, London.

[Pl2] Plotkin B. {\it Varieties of algebras and algebraic varieties.
Categories of algebraic varieties,} Siberian Advances in Mathematics, {\bf 7}
(2), , 1997

[Pl3] Plotkin B. {\it Seven lectures on the Universal Algebraic Geometry,}
Preprint \#1. Jerusalem, Inst. of Math., Hebrew University, 2000/2001

[PPT] Plotkin B., Plotkin E., Tsurkov A. {\it Geometrical equivalence of
groups,} Communications in Algebra. 27(8), 1999.

[Sha] Shakhova S.A., {\it On the lattice of quasivarieties of nilpotent
groups of class }$2${\it ,} Siberian Advances in Mathematics, {\bf 7} (3),
1997

\end{document}